\documentclass[a4paper,reqno]{amsart}
\usepackage{eucal}

\newtheorem{theorem}{Theorem}
\newtheorem{corollary}[theorem]{Corollary}
\newtheorem{definition}[theorem]{Definition}
\newtheorem{proposition}[theorem]{Proposition}
\newtheorem{remarkth}[theorem]{Remark}
\newenvironment{remark}{\begin{remarkth}\upshape}{\hfill$\diamond$\end{remarkth}}

\newcommand{\Mpi}{T^*E}
\newcommand{\Jpi}[1][1]{J^{#1}\pi}
\newcommand{\dJpi}{J^1\pi^*}
\newcommand{\h}[1]{\mathbf{h}^{(#1)}}
\newcommand{\T}[1]{\boldsymbol{T}^{(#1)}}
\newcommand{\Dt}[1]{d_{\T{#1}}}
\newcommand{\FL}{\mathcal{F}_L{}}
\newcommand{\hFL}{\widehat{\mathcal{F}}_L{}}
\newcommand{\quand}{\quad\text{and}\quad}
\newcommand{\qquand}{\qquad\text{and}\qquad}
\newcommand{\cinfty}[1]{\mathcal{C}^{\infty}(#1)}
\newcommand{\weq}[1]{\underset{#1}{\approx}}
\newcommand{\seq}[1]{\underset{#1}{\equiv}}
\newcommand{\map}[3]{#1\colon#2\to#3}
\renewcommand{\sec}[1]{\operatorname{Sec}(#1)}
\newcommand{\pd}[2]{\frac{\partial#1}{\partial #2}}
\newcommand{\Real}{\mathbb{R}}
\newcommand{\pb}[3][]{\{#2,#3\}_{#1}}
\newcommand{\vf}[1]{\mathfrak{X}(#1)}
\newcommand{\vvf}[1]{\mathfrak{X}^V(#1)}

\begin{document}

\title[Second Noether's theorem and gauge symmetries]{On second Noether's theorem and gauge symmetries in Mechanics}

\author{Jos\'e F.\ Cari\~nena}
\address{Jos\'e F.\ Cari\~nena:
Departamento de F\'{\i}sica Te\'{o}rica, Universidad de Zaragoza,
Pedro Cerbuna 12, 50009 Zaragoza, Spain}
\email{jfc@unizar.es}

\author{Joan-Andreu L\'azaro-Cam\'{\i}}
\address{Joan-Andreu L\'azaro-Cam\'{\i}:
Departamento de F\'{\i}sica Te\'{o}rica, Universidad de Zaragoza,
Pedro Cerbuna 12, 50009 Zaragoza, Spain}
\email{lazaro@unizar.es}

\author{Eduardo Mart\'{\i}nez}
\address{Eduardo Mart\'{\i}nez:
Departamento de Matem\'atica Aplicada, Universidad de Zaragoza,
Pedro Cerbuna 12, 50009 Zaragoza, Spain}
\email{emf@unizar.es}

\begin{abstract}
We review the geometric formulation of the second Noether's theorem in time-dependent mechanics. The commutation relations between the dynamics on the final constraint manifold and the infinitesimal generator of a symmetry are studied. We show an algorithm for determining a gauge symmetry which is
closely related to the process of stabilization of constraints, both in Lagrangian and Hamiltonian formalisms. The connections
between both formalisms are established by means of the
time-evolution operator.
\end{abstract}

\keywords{Constrained Lagrangian systems, constrained Hamiltonian systems, Noether's theorem, gauge symmetries}

\maketitle

\section{Introduction and notation}

The aim of this paper is to review, from a geometric point of view, the second Noether's theorem in the framework of time-dependent mechanics as developed in \cite{CFM91}, in order to analyse the search for the generating functions of infinitesimal gauge symmetries, using the tools of modern  differential geometry. This theorem was stated in the context of field theories \cite{Noe}, but we feel convenient to analyse first the results of the theorem in the simpler case of time-dependent mechanics, before proceeding to the more general case. The second Noether's theorem can be
stated in the traditional formulation as follows:

\begin{theorem}
If the action integral $S=\int L(t,q^\alpha,v^\alpha)\,dt$ is
invariant under an infinitesimal transformation involving an arbitrary 
function and its derivatives up to order $N$, then there exists an identity relation between the Lagrangian equations $\delta L_\alpha=0$ and their derivatives.  More precisely, if the action is invariant under the infinitesimal transformation  
\begin{equation}
\delta q^\alpha=\sum_{k=0}^N\varepsilon^{\left(  k\right)
}X^\alpha_{k}(t,q^\alpha,v^\alpha) \,,\qquad \varepsilon\left(  t\right)
~\text{arbitrary function}\,,\quad \varepsilon^{\left(  k\right)  }=\frac
{d^{k}\varepsilon}{dt^{k}} \,,\label{eq 1}%
\end{equation}
then there is an identity, called Noether identity:
\begin{equation}
\label{eq 2}
\sum_{k=0}^{N}(-1)^k\frac{d^k}{dt^k}(X_k^\alpha\,\delta L_\alpha)=0.
\end{equation}
\end{theorem}

It is an easy exercise to  prove the result  of this theorem writing down the variation of the action under such an infinitesimal variation, performing  an appropriate integration by parts, and using that  the function $\varepsilon$ is arbitrary. However, we want to understand it from a more intrinsic point of view.

The geometric setting of the Lagrangian formalism for field theories uses the
theory of jet bundles~\cite{Sau}. One starts with a fibre bundle $\map{\pi}{E}{\Real}$ and for any pair of integer numbers $k,l$ with $k>l\geq 0$ defines the natural projections $\map{\pi_{k,l}}{\Jpi[k]}{\Jpi[l]}$, where $\Jpi[0]=E$ and denote $\map{\pi_{k}}{\Jpi[k]}{\Real}$ the composition $\pi_{k}=\pi\circ\pi_{k,0}$. The coordinate in $\Real$ will be denoted $\left(  t\right)$, and we will refer to it as the time, and in $E$ by $(t,q^\alpha)$. Similarly, coordinates in $\Jpi[k]$ will be denoted  $\bigl(t,q^\alpha,{q^\alpha}_{(1)},...,{q^\alpha}_{(k)}\bigr)$.  In the particular case of $k=1$, $\Jpi$, usually called evolution space, its coordinates will be denoted in the more usual way $(t,q^\alpha,v^\alpha)$ and for $k=2$,  $\Jpi[2]$,  $(t,q^\alpha,v^\alpha,a^\alpha)$. When we choose a trivialization $E\simeq\Real\times Q$ then the jet bundles also trivialize $\Jpi[k]\simeq\Real\times T^kQ$, but, in general, there is no preferred trivialization.

The dynamics is defined as follows: Given a $\pi_1$-semibasic 1-form
$\boldsymbol{L}=L\,  dt$, we define the Poincar\'e-Cartan
forms,  $\Theta_{L}=dL\circ S+L\, dt\in \bigwedge^1(J^1\pi)$, with $S$ being
the vertical endomorphism,  and 
$\Omega_{L}=-d\Theta_{L}\in \bigwedge^2(J^2\pi)$ \cite{FdP}.
The dynamics is then described by the integral curves of 
 a vector field $\Gamma_{L}$ such that
$$
i_{\Gamma_{L}}\Omega_{L} =0,
\qquad
i_{\Gamma_L}dt=1,
\qquand
S\left(\Gamma_{L}\right)=0,
$$
the last two conditions being the so-called \textsc{sode} conditions. Alternatively, a \textsc{sode} vector field $\Gamma$ can be identified with a section $\map{\gamma}{\Jpi}{\Jpi[2]}$. An integral curve of~$\Gamma$ is a local section $\map{\sigma}{\Real}{E}$ of $\pi$ such that $j^2\sigma=\gamma\circ j^1\sigma$. 

When the Lagrangian is singular \cite{Nar2}, the dynamical vector field
$\Gamma_{L}$ only exists on the primary constraint manifold $S_{L}^{1}$.
Tangency requirements on such vector may lead to new constraints and, consequently, to new constraint submanifolds, until $\Gamma_{L}$ is eventually tangent to the so-called  final constraint manifold $S_{L}^{f}$. In  such a way, solutions to the dynamics only exist at points of $S_{L}^{f}$.

\section{Jet fields and vector fields along maps\label{along}}

The necessity of considering a more general concept than that of vector field, but vector fields along a map, in order to deal in a proper way with the symmetry of the dynamics  was proved in \cite{converse}.  Similarly, 
the usual concept of Noether symmetry, an infinitesimal point transformation leaving  the action invariant up to possible border terms,  can also be extended in this more  precise geometric language  in terms of the Lagrangian
density function $L$, and not only jet prolongations of projectable vector fields will be considered, but the starting object will also be a vector field along the projection $\map{\pi_{1,0}}{\Jpi}{E}$. We first establish the concept of prolongation of such an object and  introduce vector fields along the natural projections
$\map{\pi_{k+1,k}}{\Jpi[k+1]}{\Jpi[k]}$ built from $X\in\vf{\pi_{1,0}}$. What follows is a summary of what can be found, written in a slightly different way, in \cite{And} and \cite{Laz}.

\begin{definition} 
If $X$ is a $\pi$-projectable vector field along $\pi_{1,0}$, we can define, for any integer number $k$, $X^{(k)}\in\vf{\pi_{k+1,k}}$ as follows:
\begin{equation}
X_{j_{t}^{k+1}\phi}^{\left(  k\right)  }f=\frac{d}{ds}\Big\vert_{s=0}f\left(  j_{\varphi_{s}\left(  t\right)  }^{k}\left(
\pi_{k,0}\circ\Psi_{s}\circ j^{k}\phi\circ\varphi_{-s}\right)  \right)\,,
\qquad f\in\cinfty{\Jpi[k]}  \,,\label{eq 7}%
\end{equation}
where $\Psi_{s}$ is the flow of any auxiliary vector field
$\tilde{X}\in\vf{\Jpi}$ such that $\pi_{1,0\,*}\circ\tilde{X}=X$ and $\varphi_{s}$ is the flow induced by $X$ in the base manifold $M$.
\end{definition}

One can easily check that the vector field $X^{(k)}\in\vf{\pi_{k+1,k}}$ along $\pi_{k+1,k}$ is well defined and $T\pi_{k,k-1}\circ X^{(k)} = X^{(k-1)}\circ\pi_{k+1,k}$. Then, as done in \cite{Edu}, we consider the $\pi_{k+1,k}$-derivations $i_{X^{\left(  k\right)  }}$ and $d_{X^{\left(  k\right)  }}$
such that 
 $d_{X^{\left(  k\right)  }}=i_{X^{\left(  k\right)  }}\circ d+d\circ
i_{X^{\left(  k\right)  }}$.

There is also a canonical jet field $\map{\h{k}}{\Jpi[k+1]}{J^1\pi_k}$ along $\pi_{k+1,k}$ defined by
\begin{equation}
\label{eq 9}%
\h{k}\left(  j_{t}^{k+1}\phi\right)  =j_{t}%
^{1}\left(  j^{k}\phi\right),
\end{equation}
and the associated derivation  
$\map{d_{\h{k}}}{\Omega^*(\Jpi[k])}{\Omega^*(\Jpi[k+1])}$ along $\pi_{k+1,k}$ given by   $d_{\h{k}}=i_{\h{k}}\circ d-d\circ i_{\h{k}}$,
corresponds to the classical operation of \textit{total time differential} of a differential form. For instance, acting on a function $f\in\cinfty{\Jpi[k]}$,
\begin{equation}
d_{\h{k}}\left(  f\right)  =\Dt{k}(f)\,dt=\left[\pd{f}{t}+\sum_{r=0}^{k}q^\alpha_{(r+1)}\pd{f}{q^\alpha_{(r)}}
\right]  dt.
\end{equation}
Other properties involving these derivations 
can be found in \cite{And} and \cite{Laz}.

For $k=1$, with the help of the map $\h{1}$, a \textsc{sode} section $\gamma$ defines a holonomic jet field $\map{\mathcal{Y}}{\Jpi}{J^1\pi_1}$ by composition $\mathcal{Y}=\h{1}\circ\gamma$. Jet fields obtained in this way are called holonomic jet fields. In particular, the solutions to the Euler-Lagrange equations are the integral 
curves of the holonomic jet fields $\map{\mathcal{Y}_{L}}{\Jpi}{J^1\pi_1}$ such that its horizontal projector $h_{\mathcal{Y}_{L}}$ satisfies
$
i_{h_{\mathcal{Y}_{L}}}\Omega_{L}=0
$.

\section{Noether symmetries\label{noether}}

Once the main geometrical ingredients have been introduced, we can express the infinitesimal variation  of the Lagrangian density function $L$ with respect a vector field $X\in\vf{\pi_{1,0}}$ in terms of the Poincar\'{e}--Cartan $1$-form $\Theta_{L}$ and the $2$-form of Euler-Lagrange
$\delta L = d_{\h{1}}\Theta_{L} + \pi_{2,1}^*d\boldsymbol{L} \in\Omega^{2}(\Jpi[2])$ (which locally takes the form
$\delta L=\delta L_\alpha dq^\alpha\wedge dt$, where the $\delta L_\alpha=\pd{L}{q^\alpha}-\Dt{1}\left(\pd{L}{v^\alpha}\right)$ give rise to the Euler-Lagrange equations), as 
indicated in the
following proposition:

\begin{proposition}
Given a  $\pi$-projectable $X\in\vf{\pi_{1,0}}$, then
\begin{equation}
d_{X^{(1)}}\boldsymbol{L}=i_{X_{ev}\circ\pi_{2,1}}\delta L^{\vee
}+d_{\h{1}}\left(  i_{X}\Theta_{L}^{\vee}\right)
,\label{eq 10}%
\end{equation}
where $X_{ev}=\left(  X^\alpha-v^\alpha F\right) \, {\partial}/{\partial
q^\alpha}$ is the evolutionary vector field of $X=F\pd{}{t}+X^\alpha\pd{}{q^\alpha}$ (see \cite{And} for an
intrinsic definition) and $\delta L^{\vee}\in\Omega^{2}\left(  \pi
_{2,0}\right)  $ and $\Theta_{L}^{\vee}\in\Omega^{1}\left(  \pi_{1,0}\right)
$ are the forms along $\pi_{2,0}$ and $\pi_{1,0}$ associated with the semibasic
forms  $\delta L$ and
$\Theta_{L}$ respectively (see \cite{CFM92}).
\end{proposition}

The concept of Noether symmetry can be introduced as follows:

\begin{definition}
[\cite{CFM92}]We say that a $\pi$-projectable $X\in\vf{\pi_{1,0}}$
 is a Noether symmetry if there exists a function
 $F\in\cinfty{\Jpi}$ such that
\begin{equation}
d_{X^{(1)}}\boldsymbol{L}=d_{\h{1}}F\,,
\end{equation}
or in an equivalent way, if there is a function
$G\in\cinfty{\Jpi}$, such that
\begin{equation}
\label{eq 4}
i_{X_{ev}\circ\pi_{2,1}}\delta L^{\vee}=d_{\h{1}}G.
\end{equation}
\end{definition}

Observe that if $X\in\vf{\pi_{1,0}}$ is a Noether
symmetry, then $G=F-i_{X}\Theta_{L}^{\vee}$ is a conserved quantity. This
is the content of the first Noether's theorem. As a consequence of 
the 
proposition,  we can restrict the study of symmetries $X\in\vf{\pi_{1,0}}$ to those that are $\pi$-vertical, since we are
only interested in the evolutionary part of $X$, which is $\pi$-vertical. We
will indicate this fact by $X\in\vvf{\pi_{1,0}}$.

It is worth noticing that if $X\in\mathfrak{X}^{V}\left(  \pi_{1,0}\right)  $ is an exact Noether symmetry, i.e. $i_{X_{_{ev}\circ\pi_{2,1}}}\delta L^{\vee }=d_{\mathbf{h}^{\left(  1\right)  }}G,$ then $G$ is $\mathcal{F}_{L}%
$-projectable; in other words, $G=\mathcal{F}_{L}^{\ast}\left(  G_{H}\right)
$ for some $G_{H}\in\mathcal{C}^{\infty}\left(  J^{1}\pi^{\ast}\right)$ (see
\cite{Pons}).

When the Lagrangian is singular, we slightly modify 
this definition in order to take into account  the relevant r\^ole played
by the submanifold
$S_{L}^{f}$. Two forms $\omega$ and $\eta$ are said to be 
weakly equivalent, and we will write $\omega\weq{S_L^f}\eta$,  when their  pull-backs coincide on the final constraint
manifold. In the
same way,  they are said to be strongly equivalent if both the forms and their differentials are weakly equivalent. In this case, we will write
$\omega\seq{S_L^f}\eta$. So, given a $\pi$-vertical
$X\in\vf{\pi_{1,0}}$, it is said to be either (1) an exact Noether symmetry if
$i_{X}\delta L^{\vee}=d_{\h{1}}G$, (2) a weak
Noether symmetry if $i_{X}\delta L^{\vee}\weq{S_L^f}d_{\h{1}}G$, or (3) a strong Noether symmetry
when $i_{X}\delta L^{\vee} \seq{S_L^f} d_{\h{1}}G$. These distinctions will be useful later.

\section{Hamiltonian formalism}

In contrast with the autonomous case, in time-dependent Hamiltonian Mechanics the Hamiltonian is not a function but a section $h$ of a certain bundle. Given a bundle $\map{\pi}{E}{\Real}$ we consider the affine-dual bundle $\operatorname{Aff}(\Jpi,\Real)$, which is canonically isomorphic to $T^*E$, and also the vector bundle $\map{p}{\dJpi\equiv\operatorname{Ver}(\pi)^*}{E}$ dual to the vertical bundle. We have an affine bundle fibration $\map{\mu}{\Mpi}{\dJpi}$ and a Hamiltonian is a section $h$ of the projection $\mu$. The associated vector bundle is $p_1^*(T^*\Real)\to\dJpi$, where $p_1$ is the projection $\map{p_1}{\dJpi}{\Real}$. Since we have a canonical 1-form $dt$ on $\Real$, the bundle $p_1^*(T^*\Real)$ is canonically isomorphic to $\dJpi\times\Real\to\dJpi$, the isomorphism being $(\beta,a\,dt)\mapsto (\beta,a)$. Therefore, a section of the associated vector bundle will be considered as a function on $\dJpi$.

Given a Hamiltonian section $h\in\sec{\mu}$, the pullback by $h$ of the canonical symplectic form $\Omega$ on $T^*E$ defines a 2-form $\Omega_h=h^*\Omega$ on $J^1\pi^*$. The associated Hamiltonian vector fields are the solutions $\Gamma_h$ to the equations
\begin{equation}
i_{\Gamma_h}\Omega_{h}=0
\qquand
i_{\Gamma_h}dt=1,
\end{equation}
and the solution to the Hamilton equations are the integral curves of such vector fields.

The relation with the Lagrangian formalism is as follows (see \cite{Car2} for details). From the Lagrangian $L$ we can define two maps, usually called the Legendre transformation $\map{\FL}{\Jpi}{\dJpi}$ and the extended Legendre transformations $\map{\hFL}{\Jpi}{\Mpi}$, related by $\mu\circ\hFL=\FL$. In local coordinates $(t,q^\alpha,u,p_\alpha)$ in $\Mpi$ and $(t,q^\alpha,p_\alpha)$ in $\dJpi$ adapted to the fibration $\mu$ we have that 
$$
\hFL(t,q^\alpha,v^\alpha)=\left(t,q^\alpha,-E_L,\pd{L}{v^\alpha}\right)
\quand
\FL(t,q^\alpha,v^\alpha)=\left(t,q^\alpha,\pd{L}{v^\alpha}\right).
$$
In this expression $E_L$ is the energy function $E_L=v^\alpha\pd{L}{v^\alpha}-L$ defined by $L$.

When the Lagrangian is hyper-regular we have that $\FL$ is invertible and a unique section $h$ of $\mu$ is determined by the equation $\hFL=h\circ\FL$. In the singular case, any section $h$ satisfying the above relation will be called a Hamiltonian for $L$, and therefore, if it exists, $h$ needs only to be defined on the primary constraint manifold $S_H^1=\FL(\Jpi)$.

Throughout this paper we will assume that the Lagrangian $L$ is such that a Hamiltonian section $h$ exists. Two Hamiltonian sections $h$ and $h'$ differ by a section of the associated vector bundle $p_1^*(T^*\Real)$, which can be identified with a function $\phi$ on $\dJpi$. From $\hFL=h\circ\FL$ and $\hFL=h'\circ\FL$ with $h'=h+\phi$, we get $\phi\circ\FL=0$. In other words, the difference $\phi$ between two Hamiltonian sections is a primary constraint. Thus if $\{\phi_\alpha\}$ is a complete set of functionally independent primary constraints then we can write a general Hamiltonian section for $L$ in the form 
\begin{equation}
\label{eq 11}
h=h_c+\lambda^I\phi_I
\end{equation}
where $h_c$ is a particular one. As in the Lagrangian case, tangency requirements will lead to new constraints until we eventually find a solution $\Gamma_h$ tangent to the final constraint manifold $S_H^f$.

In the Hamiltonian formalism for autonomous mechanics, there are the concepts of first and second class functions (see for example \cite{Car3}). In order to extend these concepts to the time-dependent case, we will introduce a Poisson bracket in $\dJpi$. As we will see later, symmetries which are Hamiltonian vector fields of some functions will be well linked to symmetries in the Lagrangian formalism.

Let $\Omega$ be the canonical symplectic form on $\Mpi$ and denote by $\pb[\Mpi]{~}{~}$ the associated Poisson bracket. Given two functions $f$, $g\in\cinfty{\dJpi}$, it is easy to see that $\pb[\Mpi]{\mu^*f}{\mu^*g}$ is a $\mu$-projectable function, so that following \cite{Sar}, we can define a unique Poisson structure $\pb[\dJpi]{~}{~}$ on $\dJpi$ such that $\mu$ is a Poisson map. In other words, the Poisson bracket on $\dJpi$ is defined by  
\begin{equation}
\mu^*\pb[\dJpi]{f}{g}=\pb[\Mpi]{\mu^*f}{\mu^*g}.
\end{equation}
In this way, we can associate to every function $g\in\cinfty{\dJpi}$ its Hamiltonian vector field 
$$
Y_g=\pb[\dJpi]{\ \ }{g} =\pd{g}{p_\alpha}\pd{}{q^\alpha}-\pd{g}{q^\alpha}\pd{}{p_\alpha}.
$$

Given a section $h\in\sec{\mu}$, for every $\alpha\in\Mpi$, we have that the linear forms $\alpha$ and $h(\mu(\alpha))$ both project to $\mu(\alpha)$ so that they differ in an element of the vector bundle $\alpha-h(\mu(\alpha)) = h^\star(\alpha)\,dt$. In this way, we have associated an affine function $h^\star\in\cinfty{\Mpi}$ to the section $h\in\sec{\mu}$. Moreover, the canonical Poisson bracket of two functions of this type is projectable to a function on $\dJpi$. The Hamiltonian vector field $X_{h^\star}=\pb[\Mpi]{~}{h^\star}\in\vf{\Mpi}$ associated to $h^\star$ is $\mu$-projectable and projects to the vector field $X_h\in\vf{\dJpi}$, determined by 
\begin{equation}
\mu^*(\mathcal{L}_{X_h}f)= \pb[\Mpi]{\mu^*f}{h^\star}.
\end{equation}
It is easy to see that if $h'=h+g$ then $h'{}^\star=h^\star+\mu^*g$ and thus $X_{h'}=X_h+Y_g$. Also notice that $i_{X_h}\Omega_h=0$ on the primary constraint manifold, but in general $X_h$ is not tangent to such manifold.

In local coordinates, the section $h$ is of the form $h(t,x^i,p_i)=(t,x^i,-H(t,x,p),p_i)$, the function $h^\star$ is $h^\star(t,x,u,p)=u+H(t,x,p)$ and the Hamiltonian vector field $X_h$ is
\begin{equation}
X_h=\pd{}{t}+\pd{H}{p_\alpha}\pd{}{q^\alpha}-\pd{H}{q^\alpha}\pd{}{p_\alpha}.
\end{equation}
Therefore, the dynamic evolution of a function may be expressed in terms of a
Poisson bracket, in the same way as in time-independent mechanics. In local coordinates%
\begin{equation}
\dot{f}=X_hf
=\pd{f}{t}+\pd{f}{q^\alpha}\pd{H}{p_\alpha}-\pd{f}{p_\alpha}\pd{H}{q^\alpha}.
\end{equation}
which is of the classical form $\dot{f}=\pd{f}{t}+\pb[\dJpi]{f}{H}$. Notice however that the vector field $\pd{}{t}$ and the function $H$ are only locally defined on $\dJpi$.

\begin{remark}
The above properties can also be alternatively understood as follows. We have a Lie algebroid structure over the affine bundle $\map{\mu_1}{J^1\mu}{\dJpi}$. (See~\cite{MaMeSa, GGU} for the details.)  The anchor is the affine map $\map{\rho}{J^1\mu}{T\dJpi}$ defined by $\rho(j^1_\beta h)=X_h(\beta)$, for every section $h$ of $\mu$ and every $\beta\in\dJpi$ and the associated linear map  $\map{\vec{\rho}}{T^*\dJpi}{T\dJpi}$ is determined by $\rho(dg)=Y_g$, that is by the Poisson tensor on $\dJpi$. The action of a section $j^1h$ of $J^1\mu$ on a section $\phi$ of the associated vector bundle is by Lie derivative with respect to the vector field $X_h$.
\end{remark}

\begin{definition}
A function $g\in\cinfty{\dJpi}$ is said to be first class function if, for every (final) constraint function $\phi$, we have $\pb[\dJpi]{g}{\phi}\weq{S_H^f}0$. A section $\eta\in\sec{\mu}$ is said to be a first class section if $\mathcal{L}_{X_\eta}\phi\weq{S_H^f}0$ for every (final) constraint function $\phi$. A section of $\mu$ or a function of $\dJpi$ which is not of first class is said to be second class.
\end{definition}

Therefore, we can take a complete set of independent primary constraints $\{\phi_I\}$ and partition it $\{\phi_I\}=\{\phi_a,\phi_\alpha\}$ into first class constraints $\{\phi_a\}$ and second class constraints $\{\phi_\alpha\}$. When we apply the constraint algorithm, demanding that all the constraints have to be preserved by $\Gamma_{h}$, we can determine in Eq. (\ref{eq 11}) the multipliers $\lambda^\alpha$ which correspond to second class constraints. On the contrary, those multipliers $\lambda^a$ corresponding to primary first class constraints will remain indeterminate, being a sign of the non-uniqueness of the solution to the dynamics.

Once we have fixed the multipliers $\lambda^\alpha$ corresponding to second class primary constraints and we have got a Hamiltonian section $h_0=h_c+\lambda^\alpha\phi_\alpha$, we can only add first class primary constraints. Therefore any Hamiltonian for $L$ is of the form $h=h_0+\lambda^a\phi_a$ where $\{\phi_a\}$ is a complete set of functionally independent first class primary constraints.

\begin{definition}
Given a section $h\in\sec{\mu}$ such that $X_{h}$ is tangent to the primary constraint manifold $S_H^1$, a function $g\in\cinfty{\dJpi}$ is said to generate a Hamiltonian (Noether) symmetry
if $\pb[\Mpi]{\mu^*g}{h^{\star}}\weq{\tilde{S}_H^2}$, where $\tilde{S}_{H}^2=\mu^{-1}(S_{H}^2).$
\end{definition}

The reason why we introduce the concept of a Hamiltonian symmetry by means of
the vanishing of $\pb[\Mpi]{\mu^*g}{h^{\star}}$ on
$\tilde{S}_{H}^{2}$ will be clear in the following sections. The idea is that
this kind of symmetries are in correspondence with Noether symmetries in the
Lagrangian framework.

\section{The time evolution operator}

In addition to the Legendre transformation, there is a second geometric object, the time evolution operator, that connects the Lagrangian and Hamiltonian formalisms. 

\begin{theorem}
[\cite{Car4}]Given a Lagrangian $L\in\cinfty{\Jpi}$, there exists a unique vector field $\widehat{K}$ along $\hFL$, called the extended time
evolution operator, such that
\begin{enumerate}
\item $T\pi_{E}\circ\widehat{K}=i_1$, where $\map{i_1}{\Jpi}{TE}$ is the map $i_1(j^1_t\sigma)=\dot{\sigma}(t)$, and
\item $i_{\widehat{K}}\Omega=0$, where  $\Omega$ is the canonical symplectic form on $T^*E$.
\end{enumerate}
\end{theorem}

In local coordinates,%
\begin{equation}
\widehat{K}=\left(\pd{}{t}\circ\hFL\right)
+
v^i\left(\pd{}{q^\alpha}\circ\hFL\right) +\pd{L}{t}\left(\pd{}{u}\circ\hFL\right)
+\pd{L}{q^\alpha}\left(\pd{}{p_\alpha}\circ\hFL\right).
\end{equation}
By composition with the differential of the projection $\mu$ we get a vector field $K=T\mu\circ\widehat{K}$ along $\FL$, the (restricted) time evolution operator, whose coordinate expression is 
\begin{equation}
K=\left(\pd{}{t}\circ\FL\right)
+
v^i\left(\pd{}{q^\alpha}\circ\FL\right) 
+\pd{L}{q^\alpha}\left(\pd{}{p_\alpha}\circ\FL\right).
\end{equation}

\begin{remark}[Alternative construction of $\hat{K}$]
Let $L\in\cinfty{\Jpi}$ be a Lagrangian and consider the associated homogeneous Lagrangian $\hat{L}\in\cinfty{\mathaccent23{T}E}$ given in coordinates by 
$
\hat{L}(x^0,x^i,w^0,w^i)=w^0\,L(x^0,x^i,w^i/w^0)
$.
We consider the autonomous time evolution operator $K_{\hat{L}}$ defined by $\hat{L}$, that is $K_{\hat{L}}=\chi\circ d\hat{L}$ where $\chi$ is the canonical isomorphism $\map{\chi}{T^*(TE)}{T(T^*E)}$. By composition with the canonical inclusion $\map{i_1}{\Jpi}{TE}$ we get the vector field $\hat{K}=K_{\hat{L}}$ along $\hFL$. The above properties are easy to prove from this definition.
\end{remark}

The integral curves of $K$ are the solutions to the Euler-Lagrange equations for $L$. This is clear from the following relation between the time evolution operator and the total time derivative.

\begin{proposition}
\label{prop1}
Let $G\in\cinfty{\dJpi}$ and denote by $X\in\vvf{\pi_{1,0}}$ the vector field given by $X(v)=Tp\,(Y_G(\FL(v)))$, for $v\in\Jpi$. Then 
\begin{equation}
\label{eq 3}
d_{\h{1}}(\FL^*G)=\pi_{2,1}^*(K\cdot G)\,\pi_{2}^*(dt) -i_{X}\delta L^{\vee}.
\end{equation}
\end{proposition}

Notice that the action of $X$ on functions $f\in\cinfty{E}$ is given by 
\begin{equation}
\label{eq 13}
X\left(f\right)=\FL^*\pb[\dJpi]{p^*f}{G} 
\end{equation}
and locally $X=\FL^*\left(\pd{G}{p_\alpha}\right)\pd{}{q^\alpha}$.

When we compose in the above expression with any section $\gamma$ solution to the dynamics we immediately get the following consequence.

\begin{proposition}
\label{prop2}
For any admissible dynamic vector field $\Gamma_{L}$
defined on the primary constraint manifold $S_{L}^{1}$ we have
\begin{equation}
\Gamma_{L}\left(\FL^*G\right)\weq{S_{L}^{1}}K\cdot G,
\end{equation}
for every function $G\in\cinfty{\dJpi}$.
\end{proposition}

In particular, since the pullback by $\FL$ of a Hamiltonian constraint is a Lagrangian constraint we get:

\begin{corollary}
If $\phi$ is a Hamiltonian constraint then $K\cdot\phi$ is a Lagrangian constraint.
\end{corollary}

In the light of the above proposition~\ref{prop1}, and taking into account our definition of the different types of Noether symmetry,  we will say that a function $G\in\cinfty{\dJpi}$ generates  
(1) an exact Noether symmetry if $K\cdot G=0$, 
(2) a weak Noether symmetry if $K\cdot G\weq{S_L^f}0$,  and
(3) a strong Noether symmetry if $K\cdot G\seq{S_L^f}0$.

Finally, we state another proposition that relates the action of $K$ on a
function $G\in\cinfty{\dJpi}$ to its dynamic evolution with respect to the part of the Hamiltonian $h_0$ determined at the first stage, and to some other terms involving primary constraints, on one hand those being first class and, on the other, those being second class. The steps taken in \cite{Pons} in this direction, in the time-independent case, have served as a guideline.

\begin{proposition}
\label{prop3}
Let $\{\phi_I\}=\{\phi_a,\phi_\alpha\}$ be a complete set of independent primary constraints with $\{\phi_a\}$ first class and $\{\phi_\alpha\}$ second class (at primary level). Then, for every function $G\in\cinfty{\dJpi}$,
\begin{equation}
\label{eq 12}
K\cdot G  =\hFL^*\pb[\Mpi]{\mu^*G}{h_0^\star}
  +\chi_\alpha M^{\alpha\beta}\FL^*\pb[\dJpi]{G}{\phi_\beta}  
  +\nu^{a}\FL^*\pb[\dJpi]{G}{\phi_a},
\end{equation}
where $M^{\alpha\beta}$
is the inverse matrix of $\FL^*\pb[\dJpi]{\phi_\alpha}{\phi_\beta}$,  $\chi_\alpha=K\cdot\phi_\alpha$ are primary Lagrangian constraints and
$\nu^a\in\cinfty{\Jpi}$ are non-projectable.
\end{proposition}

This proposition will be used later to establish the connection between gauge
symmetries in the Hamiltonian formalism and the Lagrangian one. However, it
has an immediate but important corollary, extending the results of
\cite{Pons2}.

\begin{corollary}
At primary level, a Hamiltonian constraint $\phi$ is a first class constraint  if and only if $K\cdot\phi$ is a projectable Lagrangian constraint.
\end{corollary}

\begin{remark}
We can recast proposition~\ref{prop3}, and some steps in its proof, in the following way. Given a Hamiltonian section $h$ for $L$ and a complete set of primary constraints $\{\phi_I\}$ there exist functions $\nu^I$ such that 
\begin{equation}
\label{eq 5}
K=(X_h\circ\FL)+\nu^I(Y_{\Phi_I}\circ\FL).
\end{equation}
If $\Gamma_I\in\vf{\Jpi}$ are the vector fields $\Gamma_I=(Tp\circ Y_{\Phi_I}\circ\FL)^V$ (which generate the kernel of $T\FL$), then the functions $\nu^I$ satisfy $\mathcal{L}_{\Gamma_I}\nu^J=\delta^J_I$, and hence they are non $\FL$-projectable.

If we partition $\{\Phi_I\}=\{\Phi_a,\Phi_\alpha\}$ into first class constraints $\{\phi_a\}$ and second class constraints $\{\phi_\alpha\}$, then tangency considerations imply that the multipliers $\nu^\alpha$ corresponding to the second class constraint (at primary level) are determined by $\nu^\alpha=M^{\alpha\beta}[K\cdot\phi_\beta-\FL^*(X_h\phi_\beta)]$. Moreover, by modifying the Hamiltonian section $h\to h-\lambda^\alpha\Phi_\alpha$ with $\lambda^\alpha=M^{\alpha\beta}\FL^*(X_h\phi_\beta)$ then we have the simpler form $\nu^\alpha=M^{\alpha\beta}(K\cdot\phi_\beta)$.
%
\end{remark}

\section{Commutation relations}
\label{commutation}

In this section we are going to investigate the Lie brackets of the dynamical
vector fields (in both the Lagrangian and the Hamiltonian formalisms) and the vector field that generates the symmetry. We will assume that we already performed the constraint algorithm and we have arrived to a consistent solution. 

In both, the Lagrangian and the Hamiltonian formalisms, this algorithm determines a final constraint submanifold and the set of solutions is an affine family of vector fields, that is, are sections of an affine subbundle of the tangent bundle. It follows that the concept of symmetry that we have to use is that of a symmetry of an affine subbundle. 

In all generality, we consider a manifold $N$ and an affine subbundle $\mathcal{A}\subset TN$ modeled on the vector bundle $\mathcal{V}\subset TN$. By an infinitesimal symmetry of $\mathcal{A}$ we mean a vector field $Y\in\vf{N}$ such that $[Y,\Gamma]\in\sec{\mathcal{V}}$ for every section $\Gamma\in\sec{\mathcal{A}}$.  If we fix a section $\Gamma_0\in\sec{\mathcal{A}}$ and a basis $\{\Gamma_\mu\}$ of sections of the underlying vector subbundle $\mathcal{V}$, then a section of $\mathcal{A}$ can be written in the form $\Gamma=\Gamma_0+\alpha^\mu\Gamma_\mu$. A vector field $X$ is a symmetry of $\mathcal{A}$ if and only if for every set of functions $\alpha^\mu$ there exist functions $\beta^\mu$ such that $[X,\Gamma_0+\alpha^\mu\Gamma_\mu]=\beta^\mu\Gamma_\mu$.

In the cases we are interested in, the manifold $N$ is a submanifold of a given manifold $M$. In such cases, instead of vector fields on $N$ we take vector fields on $M$ which are tangent to $N$. With this in mind, by an infinitesimal symmetry of $\mathcal{A}\subset TN\subset TM$ we mean a vector field $X\in\vf{M}$ tangent to $N$ and such that $[Y,\Gamma]\big|_N\in\sec{\mathcal{V}}$ for every vector field $\Gamma\in\vf{M}$ such that $\Gamma|_N\in\sec{\mathcal{A}}$. 

More concretely, in the cases we have in mind, the manifold $N$ is the final constraint manifold and the subbundle $\mathcal{A}$ is the set of all vectors tangent to curves solution to the dynamics on such final constraint manifold. In the Lagrangian formalism, the vector fields solution to the dynamics on the final constraint manifold are of the form $\Gamma=\Gamma_0+\alpha^\mu\Gamma_\mu$, where $\Gamma_0$ is a given \textsc{sode} and $\Gamma_\mu$ are vertical vector fields on $\Jpi$. A dynamical symmetry is a vector field $X\in\vf{\Jpi}$ tangent to the final constraint manifold $S^f_L$ such that $[X,\Gamma]\weq{S^f_L}\beta^\mu\Gamma_\mu$ for some functions $\beta^\mu\in\cinfty{\Jpi}$.
On the Hamiltonian counterpart, the vector fields solution to the dynamics are of the form $\Gamma=X_{h_0}+\alpha^\mu X_{\phi_\mu}$, where $h_0$ is a given solution and $\{\phi_\mu\}$ is a complete set of primary first class constraints (at the final level). A symmetry of the dynamics is a vector field $X\in\vf{\dJpi}$ such that $[X,\Gamma]\weq{S^f_H}\beta^\mu X_{\phi_\mu}$ for some functions $\beta^\mu\in\cinfty{\dJpi}$.

\medskip

From now on, by a first class (respectively, second class) function we mean a function which is first class (respectively, second class) with respect to the final set of constraints, i.e. with respect to the final constraint manifold manifold. Also, in what follows in this paper we will denote by $\tilde{S}^k_H\subset\Mpi$ the submanifold $\tilde{S}^k_H=\mu^{-1}(S^k_H)$.

\medskip

We consider first the Lagrangian formalism. Since the symmetry $X\in\vf{\pi_{1,0}}$ is a vector field along $\pi_{1,0}$ and any solution $\Gamma_{L}\in\vf{\Jpi}$, is a vector field on $\Jpi$, we need to prolong $X$ to a vector field defined on $\Jpi$, as well. As we previously said, a \textsc{sode} can be identified with a section $\map{\gamma}{\Jpi}{\Jpi[2]}$ of $\pi_{2,1}$. If $L$ is singular, we may have several sections $\gamma$ solution to the dynamical equation defined on the primary constraint submanifold $S_{L}^{1}$.

Given a vector field $U\in\vf{p}$ along the projection $p$, we consider the vector field $X=U\circ\FL\in\vf{\pi_{1,0}}$ along the projection $\pi_{1,0}$. For any \textsc{sode} section $\gamma$ we may consider the vector field $X^{(1)}\circ\gamma\in\vf{\Jpi}$. In particular, for solutions $\gamma$ of the dynamics, from proposition~\ref{prop2} we have:

\begin{proposition}
The vector field $X^{(1)}\circ\gamma$ does not depend on the particular choice of the section $\gamma$ solution to the Lagrangian dynamical equations on $S^1_L$. In coordinates, if $U=U^\alpha(t,q^\alpha,p_\alpha)\pd{}{q_\alpha}$ then 
$$
X^{(1)}\circ\gamma
=\FL^*(U^\alpha)\pd{}{q^\alpha}+(K\cdot U^\alpha)\pd{}{v^\alpha}.
$$
\end{proposition}

In particular, given a function $G\in\cinfty{\dJpi}$ generating a (exact, strong or weak) Noether symmetry , we consider $U_G=Tp\circ Y_G$, where $Y_G$ is the Hamiltonian vector field defined by $G$, and we apply the above procedure, arriving to a vector field $Z_G=(Tp\circ Y_G\circ\FL)^{(1)}\circ\gamma$. In local coordinates, we find
\begin{equation}
Z_G=\FL^*(\pb[\dJpi]{q^\alpha}{G})\pd{}{q^\alpha}+(K\cdot\pb[\dJpi]{q^\alpha}{G})\pd{}{v^\alpha}.
\end{equation}
In \cite{Gra} the construction of this vector field was carried out in a different but equivalent way. Now, we can establish

\begin{theorem}
\label{theo1}
Consider the following different conditions for a function
$G$ on~$\dJpi$:
$$
  (i)~K\cdot G\seq{S_H^f}0,\qquad
 (ii)~K\cdot G\weq{S_H^1}0,\qquad
(iii)~K\cdot G=0.
$$
Then,
\begin{enumerate}
\item The vector field $X\in\vvf{\pi_{1,0}}$ is a Noether symmetry with
conserved quantity $G_{L}$ if and only if the function $G$ such that $\FL^*G=G_L$ satisfies condition~$(ii)$.

\item Condition $(ii)$ holds if and only if $G$ generates a
Hamiltonian (Noether) symmetry.

\item If condition $(iii)$ holds then $\FL{}_*(Z_{G})=Y_G$.

\item If condition $(ii)$ holds then $\FL{}_*(Z_{G})\weq{S^f_H}Y_G$.

\item If $G$ is first class and condition $(i)$ holds then  $Z_G$ is a dynamic symmetry, i.e. $[Z_{G},\Gamma]\weq{S_L^f}\beta^\mu\Gamma_\mu$, for every final solution $\Gamma$.
\end{enumerate}
\end{theorem}

\begin{proof}
[Sketch of the proof]%
\item[ (1)] 
$(\Leftarrow)$ 
The primary Lagrangian constraints may be written as $\chi_{I}dt=i_{Z_{\phi_{I}}}\delta L$, so that there exist $r_I\in\cinfty{\Jpi}$ such that $K\cdot G=r^{I}\chi_{I}$. It follows from Eq.
(\ref{eq 3}) that  $X=Z_{G}-r^{I}Z_{\phi_{I}}\in\vvf{\pi_{1,0}}$ is a Noether symmetry.

\noindent  $\Rightarrow)$ From equations (\ref{eq 4}) and (\ref{eq 3}) we deduce that if $X\in\vvf{\pi_{1,0}}$ is a Noether symmetry then $\pi_{2,1}^*(K\cdot G) dt=i_{(Z_G-X)}\delta L^{\vee}$. Since the left hand side of this relation does not depend on accelerations $i_{(X-Z_G)}\delta L^{\vee}$ is a primary Lagrangian constraint.

\item[ (2)] 
$\Leftarrow)$ See the characterization given by the second item of the theorem below. The result follows by taking into account  Eq. (\ref{eq 12}) and that the pullback by $\FL$ of a secondary Hamiltonian constraint is a primary Lagrangian constraint.

\noindent 
$\Rightarrow)$ Consider a complete set of independent primary first class constraint $\{\phi_{a}\}$ and the corresponding vector fields $\Gamma_{a}\in\ker(T\FL)$. Writing condition $(ii)$ in the form $K\cdot G=r^I\chi_I$ and applying $\Gamma_a$ to it we obtain
$\FL^*\pb[\dJpi]{G}{\phi_{a}} =(\Gamma_{a} r^{I})\chi_{I}$, where we have taken into account  Eq. (\ref{eq 5}). Since the left
hand side is $\FL$-projectable so does the right hand side. But $\chi_I$ are primary Lagrangian constraints, so that only $\FL$-projectable primary constraints can appear. In other words, only the pull-back by $\FL$ of secondary Hamiltonian constraints can appear (see \cite{Pons}), and thus $\pb[\dJpi]{G}{\phi_a}\weq{S_H^2}0$. The result follows by taking into account Eq. (\ref{eq 12}), i.e. $\hat\FL^*\pb[\Mpi]{\mu^*G}{h_0^\star}\weq{S_L^1}0$, which implies $\pb[\Mpi]{\mu^*G}{h_0^\star}\weq{\tilde{S}_H^2}0$, and item $(2)$ of
theorem~\ref{theo2}.

\item[ (4)] 
We just have to show that $T\FL{}\circ Z_G$ and $Y_G\circ\FL$ agree on $S_H^f$ when acting on the variables $p_i$. For every $v\in S^f_L\subset\Jpi$,
\begin{align*}
T\FL(Z_G(v))  p_{i}
&=Z_G(v)\left(\pd{L}{v^i}\right)
=\left(\pd{}{v^i}(K\cdot G)-\FL^*\left(\pd{G}{q^i}\right)\right)(v)\\
&\weq{S_L^f}(\FL^*\pb{p_i}{G})(v)  
=\pb{p_i}{G}(\FL(v))
=Y_G(\FL(v))p_i.
\end{align*}

\item[ (3)] 
The proof is similar to that of $(4)$ and will be omitted.

\item[ (5)] 
The proof is based on the relation $Z_G(K\cdot F)=K\cdot(\pb{F}{G})+Z_F(K\cdot G)$,  for $G,F\in\cinfty{\dJpi}$. If $K\cdot G\seq{S_H^f}0$, $Z_G(K\cdot F)\weq{S_H^f} K\cdot (Y_GF)$, and if $G$ is first class, $Z_G$ is tangent to $S_L^f$. If $\{\psi_I\}$ is a complete set of Hamiltonian constraints then $\{K\cdot\psi_I\}$ is a complete set of Lagrangian constraints. Then
\begin{equation}
Z_{G}\left(  K\cdot\psi_{I}\right)  \underset{S_{L}^{f}}{\approx}K\cdot\left(
Y_{G}\left(  \psi_{I}\right)  \right)  \underset{S_{L}^{f}}{\approx}%
K\cdot\left(  \text{Ham. constraints}\right)  \underset{S_{L}^{f}}{\approx}0.
\end{equation}
On the other hand, if $F\in\cinfty{\dJpi}$, using proposition~\ref{prop2},
\begin{equation}%
\begin{aligned}
Z_G(\Gamma_L(\FL^*F))
&\weq{S_L^1}Z_G(K\cdot F)
\weq{S_L^f}K\cdot (Y_GF)\\
&\weq{S_L^f}\Gamma_L(\FL^*(Y_G F))
\weq{S_L^f}\Gamma_L(Z_G(\FL^*F)),
\end{aligned}
\end{equation}
so the Lie bracket $\left[  Z_{G},\Gamma_{L}\right]  $ is zero acting on
$\mathcal{F}_{L}$-projectable functions on the final constraint manifold.
Thus $\left[  Z_{G},\Gamma_{L}\right]$ is tangent to $S_H^f$ and it is in $\ker(T\FL)$, which concludes the proof. 
\end{proof}

Sometimes it is useful to treat separately that part of the Hamiltonian section $h$ that is determined $h_{0}$ and the primary first class term.

\begin{theorem}
\label{theo2}
Let $G\in\cinfty{\dJpi}$ be a first class function and let $\mathbf{pfcc}$ means a primary constraint that is first class with respect to $\tilde{S}_H^f$, and $\mathbf{pfcc}^{1}$ another one
that is first class with respect to $\tilde{S}_H^1$. Then,
\begin{enumerate}
\item if for every Hamiltonian section $h$ we have that $\pb[\Mpi]{\mu^*G}{h^\star}\seq{\tilde{S}_H^f}\mathbf{pfcc}$
or, equivalently, $\pb[\Mpi]{\mu^*G}{h_0^\star}\seq{\tilde{S}_H^f}\mathbf{pfcc}$ and $\pb[\Mpi]{\mu^*G}{\mathbf{pfcc}}\seq{\tilde{S}_H^f}\mathbf{pfcc}$,  with $h_0$ such that $X_{h_0}$ is tangent to $S_H^f$, then $Y_G$ is a dynamical symmetry.

\item $\pb[\Mpi]{\mu^*G}{h_0^\star}\weq{\tilde{S}_H^2}0$ and $\pb[\Mpi]{\mu^*G}{\mathbf{pfcc}^{1}}\weq{\tilde{S}_H^2}0$, with $h_0$ such that $X_{h_0}$ is tangent to $\tilde{S}_H^1$ if and only if $G$ generates a Hamiltonian (Noether) symmetry.
\end{enumerate}
\end{theorem}

\begin{proof}
\item[ (1)]
Since the difference between any section $h$ defined on $S_H^f$ and
$h_0$ (with all the multipliers associated to primary second class constraints being determined) is a primary first class constraint (with respect to $S_H^f$), the two conditions
in the statement are equivalent to $\pb[\Mpi]{\mu^*G}{h^\star}\seq{\tilde{S}_H^f}\mathbf{pfcc}$. On the other hand $i_{[Y_G,\Gamma_h]}\Omega =d(\pb[\Mpi]{\mu^*G}{h^\star}) \weq{\tilde{S}_H^f}d\left(\mathbf{pfcc}\right)$ and $\Omega$ is symplectic, from where we deduce that $[Y_G,\Gamma_h]\weq{S_H^f}\beta^a Y_{\phi_a}$, where $\{\phi_a\}$ is a set of primary first class constraints. In other words, $Y_{G}$ is a dynamic symmetry. Finally, using Eq. (\ref{eq 12}), we get $K\cdot G\seq{S_L^f}0$.

\item[ (2)]
On the primary constraint manifold the difference between $h$ and
$h_0$ (with all the multipliers associated to primary second class constraints with respect to $S_H^1$ determined) is a primary first class constraint (with respect to $S_H^1$). Thus,
$G$ generates a Hamiltonian (Noether) symmetry, $\pb[\Mpi]{\mu^*G}{h^\star}\weq{\tilde{S}_H^2}0$, if and only if  $\pb[\Mpi]{\mu^*G}{h_0^\star}\weq {\tilde{S}_H^2}0$ and $\pb[\Mpi]{\mu^*G}{\mathbf{pfcc}^1}\weq{\tilde{S}_H^2}0$.
\end{proof}

\section{Gauge symmetries}
\label{gauge}

\subsection{Lagrangian case}

In this section we will analise the existence of gauge symmetries, that is, of Noether symmetries $X\in\vvf{\pi_{1,0}}$ depending on an arbitrary function $\varepsilon(t)$ and its derivatives up to some order $N$, i.e. 
$X=\sum_{k=0}^N\varepsilon^{(k)}X_{k}$, where $X_{k}$ are vertical vector fields along $\pi_{1,0}$ and we have written $\varepsilon^{(k)}$ for the $k$-th derivative of the arbitrary function $\varepsilon(t)$. The associated conserved quantity $G_{L}$ also depends on $\varepsilon$ and its derivatives $G_{L}=\sum_{k=0}^N\varepsilon^{(k)}G_{k}^{L}$, with $G_{k}^L\in\cinfty{\Jpi}$. 

Let us suppose first that $X$ is an exact Noether symmetry. Following \cite{CFM91}, from $i_{X}\delta L^{\vee}=d_{\h{1}}G$ and taking into account that $\varepsilon$ is arbitrary, we deduce that 
\begin{equation}
\label{eq 15}
\begin{aligned}
G_{N}^{L}  & =0 \\
d_{\h{1}}G_{k}^{L}+\pi_{2,1}^*\left(
G_{k-1}^{L}\right)  dt-i_{X_{k}}\delta L^{\vee}  & =0\\
d_{\h{1}}G_{0}^{L}-i_{X_{0}}\delta L^{\vee}  &=0.
\end{aligned}
\end{equation}
for $k=1,\ldots,N$. Cari\~{n}ena et al.~\cite{CFM91} use the above relations relations~(\ref{eq 15}) to draw an algorithm for determining such symmetries. The idea is as follows. From $G_{N}^{L}=0$ and the second equation for $k=N$, we get that $i_{X_{N}}\delta L^{\vee}=\pi_{2,1}^*G^L_{N-1}$ is the pullback of a function in $\Jpi$. Thus, we must choose $X_{N}$ such that its vertical lift is in the kernel of $T\FL$ (and hence in the $\pi_{1,0}$-vertical part of the $\ker\Omega_{L}$) in order to annihilate the terms of $i_{X_{N}}\delta L^{\vee}$ which depends on the accelerations, and it follows that $G_{N-1}^{L}$ is a primary first class constraint. In this way we recover the well known fact that only singular Lagrangians may have gauge symmetries. The algorithm proceeds by choosing at every step a vector field $X_{k}\in\vvf{\pi_{1,0}}$ such that
$d_{\h{1}}G_{k}^{L}-i_{X_{k}}\delta L$ is a
$\pi_{2,1}$-projectable function, that we take as $G_{k-1}^{L}$, and it finishes when we can choose $X_{0}$ such that the difference $d_{\h{1}}G_{k}^{L}-i_{X_{0}}\delta L$ vanishes, that is $X_0$ is a Noether symmetry of $L$. 

But it may happen that, at some stage, there is no vector field $X_{k}$ such that the difference $d_{\h{1}}G_{k}^{L}-i_{X_{k}}\delta L^{\vee}$ is a $\pi_{2,1}$-projectable function. Then the algorithm fails to determine $X$ and we must restart it with a different choice of the initial vector field $X_{N}$ until we finish off all the primary constraints (in a preselected complete set of independent constraints). Thus, this algorithm has a clear problem: we do not know when it might finish successfully and, we do not know \textit{a priori} which Lagrangian primary constraints serve as input for our algorithm. 

We can reinterpret Eq. (\ref{eq 15}) from the dynamic point of view. It
reproduces the process of stabilization of a primary $\FL$-projectable constraint, $G_{N-1}^{L}$. To see this more clearly, we remark
that the conserved quantity $G_{L}$ for a Noether symmetry is always
$\mathcal{F}_{L}$-projectable, $G_{L}=\mathcal{F}_{L}^{\ast}\left(  G\right)
$. If we compose the expression $d_{\mathbf{h}^{\left(  1\right)  }}G_{k}%
^{L}+\pi_{2,1}^{\ast}\left(  G_{k-1}^{L}\right)  dt-i_{X_{k}}\delta L^{\vee
}=0$ with any admissible dynamic section $\gamma:S_{L}^{1}\subset
\Jpi\rightarrow\Jpi[2]$, it is reexpressed as $\Gamma_{L}\left(
\mathcal{F}_{L}^{\ast}\left(  G_{k}\right)  \right)  =\mathcal{F}_{L}^{\ast
}\left(  G_{k-1}\right)  ,$ so we take as $G_{k-1}^{L}$ the dynamic evolution
of the previous constraint $G_{k}^{L}$. If at some stage, $G_{k-1}^{L}$
becomes non $\mathcal{F}_{L}$-projectable, the algorithm stops. If we assume
that constraints obtained from the dynamic evolution of previous ones serve us
to characterize the next constraint manifold in the constraint algorithm, that
is to say, if no ineffective constraints appear, a constraint that never
turns into non $\mathcal{F}_{L}$-projectable is related to the indeterminacy
of the dynamics. If $\Gamma_{L}=\Gamma_{0}+\alpha^{\mu}\Gamma_{\mu}$ is the
dynamic vector field (as we denoted throughout this paper), since $\Gamma_{\mu
}\in\ker T\mathcal{F}_{L}$, if $\chi_{\nu_{0}}=\mathcal{F}_{L}^{\ast}\left(
\phi_{\nu}\right)  $ is a projectable constraint, imposing its dynamic
evolution to be zero leads to $0=\Gamma_{L}\left(  \chi_{\nu_{0}}\right)
=\Gamma_{0}\left(  \chi_{\nu_{0}}\right)  ,$ that tells us nothing about
$\left\{  \alpha^{\mu}\right\}  $. In conclusion, if there exists a gauge
symmetry, there is a primary $\mathcal{F}_{L}$-projectable constraint that never
turns into a non $\mathcal{F}_{L}$-projectable and at least one of the $\left\{
\alpha^{\mu}\right\}  $ can not be fixed (remember we have as many $\left\{
\alpha^{\mu}\right\}  $ as primary Lagrangian constraints) and $\Gamma_{L}$ is
not completely determined on $S_{L}^{f}$. For more information about the role
of projectable and non-projectable constraints see \cite{Gom} and \cite{Nar2}.

Moreover, if we chain the consecutive steps of the Lagrangian algorithm
(\ref{eq 15}) and express $G_{k}^{L}$ in terms of the previous ones, at the
end we can write%
\begin{equation}
\sum_{k=0}^{N}(-1)^k\frac{d^k}{dt^k}(X_k^\alpha\,\delta L_\alpha)=0,
\end{equation}
which is know as a Noether identity.

\subsection{Hamiltonian case}

Since the conserved quantity $G_{L}$ for a Noether symmetry is always
$\mathcal{F}_{L}$-projectable, we can look for the consequences in the
Hamiltonian formalism of having a gauge symmetry. Regarding Theorem
(\ref{theo1}) item $\left(  1\right)  $ and $\left(  2\right)  ,$ and the
second item of Theorem (\ref{theo2}), considering a Hamiltonian symmetry of
the type $G=\sum_{k=0}^{N}\varepsilon^{\left(  k\right)  }G_{k}\in
\mathcal{C}^{\infty}(J^{1}\pi^{\ast})$, where $\varepsilon^{\left(  k\right)
}$ depends only on time, implies%
\begin{equation}%
\begin{array}
[c]{rc}%
\mu^{\ast}G_{N} & \underset{\widetilde{S}_{H}^{2}}{\approx}0\\
\mu^{\ast}G_{k-1}+\left\{  \mu^{\ast}G_{k},h_{0}^{\star}\right\}  _{T^{\ast}E}
& \underset{\widetilde{S}_{H}^{2}}{\approx}0\\
\left\{  \mu^{\ast}G_{0},h_{0}^{\star}\right\}  _{T^{\ast}E} & \underset
{\widetilde{S}_{H}^{2}}{\approx}0,%
\end{array}
\label{eq 17}%
\end{equation}%
\begin{equation}
\qquad\quad\ \left\{  \mu^{\ast}G_{k},\mathbf{pfcc}^{1}\right\}  _{T^{\ast}E}%
\underset{\widetilde{S}_{H}^{2}}{\approx}0, \label{eq 18}%
\end{equation}
for $k=1,\ldots,N$. These relations suggest us a way to find the gauge
Hamiltonian (Noether) symmetries of the Hamiltonian system by means of the
following algorithm. As input, we choose for $G_{N}$ a constraint (with
respect to $S_{H}^{2}$). Then, we calculate its dynamic evolution
$\{{\mu^{\ast}G_{N}},{h_{0}^{\star}}\}_{T^{\ast}E}$ which we take as $G_{N-1}$
once all the terms involving constraints (at secondary level) have been
removed, and so on. Thus, until at some stage $G_{k}$ is nothing but a
function vanishing on $S_{H}^{2}$. If at some stage $G_{k}$ Eq. (\ref{eq 18})
is not satisfied, the algorithm can not go on and we have to try with a
different $G_{N}$.

\begin{proposition}
The operator $\hat{K}$ transforms Eqs. (\ref{eq 17}) and (\ref{eq 18}) into
Eq. (\ref{eq 15}).
\end{proposition}

A different starting point to study gauge symmetries in the Hamiltonian
formalism is adopted by Gomis et. al. in \cite{Gom} (in the context of
time-independent mechanics). They look for a dynamic gauge symmetries and,
therefore, for a $G=\sum_{k=0}^{N}\varepsilon^{\left(  k\right)  }G_{k}%
\in\mathcal{C}^{\infty}(J^{1}\pi^{\ast})$ satisfying the first item of Theorem
(\ref{theo2}). This leads them to equations equivalent to Eqs. (\ref{eq 17})
and (\ref{eq 18}) where the weak equalities $\underset{\widetilde{S}_{H}^{2}%
}{\approx}0$ must be replaced by $\underset{\tilde{S}_{H}^{f}}{\equiv
}\mathbf{pfcc}$. Essentially, under some regularity assumptions, i.e., that
the rank of the Legendre map and the matrix of Poisson brackets of constraints
is constant, and that no ineffective constraints appear, they show that is
possible to build dynamic gauge symmetries, where the functions $G_{k}$ are
first class constraints plus a quadratic term in the rest of constraints. One
can attempt to rewrite their reasonings just replacing the usual Poisson
bracket of time-independent mechanics by $\left\{  ~,\right\}  _{T^{\ast}E}$,
and everything they hold carries on being valid. Nevertheless, we emphasize
that the symmetries hence obtained are dynamic and not Noether. Actually, one
can easily verify that their symmetries satisfy $K\cdot G\underset{S_{L}^{f}%
}{\equiv}0,$ so they are not, in principle, Noether (item $\left(  2\right)  $
Theorem (\ref{theo1})). Thus, as far as we know, the problem of the existence
of gauge Noether symmetries is still unsolved.

\section{Conclusions}

We have reviewed the concept of Noether symmetry and we have clarified the
relationship between Hamiltonian constraints and Lagrangian constraints in
time-dependent mechanics, extending the results given in \cite{Pons2} for the
autonomous case. In particular, by making use of the properties of the time
evolution operator $K$, we have shown that Hamiltonian first class constraints
with respect to $S_{H}^{1}$ are transformed into $\mathcal{F}_{L}{}%
$-projectable Lagrangian constraints, and second class constraints into
non-projectable ones.

We also studied Noether symmetries that come from a function defined on the
dual affine bundle, imposing conditions on them to guarantee, in some sense,
in which cases the symmetry commutes with the dynamic vector field.
Essentially, it was done in \cite{Gra2}, but here we use, in the Lagrangian
case, the geometric content of vector fields along $\pi_{k+1,k}$ and, in the
Hamiltonian case, the new tools of Hamiltonian dynamics in the cotangent
bundle $T^{\ast}E$.

Finally, we have studied Noether symmetries depending on arbitrary functions
of the time and its derivatives. We saw that this kind of symmetries, in the
Lagrangian case, are due to the fact that there are primary $\mathcal{F}_{L}%
{}$-projectable constraints that, when imposing they have to be preserved by
the dynamic vector field, never turn into non projectable ones, being in this
way a sign of the indeterminacy of $\Gamma_{L}$. In addition, in the
Hamiltonian case, we saw we can generate Noether symmetries by means of
constraints whose dynamic evolution is always first class with respect to
primary first class at first level when regarded on $S_{H}^{2}.$ And that,
under some regularity conditions, we pointed out a way to build dynamical
symmetries that also depends on a arbitrary function and its derivatives,
following \cite{Gom}. Anyway, the existence of (Noether) gauge symmetries from
a singular Lagrangian is a problem that it is not completely solved yet.

\section*{Acknowledgments}
J. A.\ L\'{a}zaro-Cam\'{\i} acknowledges financial support from MEC (Spain) grant BES-2004-4914. Partial financial support from MEC grant BFM2003-02532 and Gobierno de Aragon grant \textsc{dga-grupos consolidados 225-206} is acknowledged.

\end{document}